\documentstyle[12pt,amssymb]{amsart}
\newcommand{\new}{\newcommand}

\new{\dt}{,\dots,}
\new{\si}{\sigma}
\new{\Si}{\Sigma}
\new{\h}[1]{\widehat {#1}}
\new{\ol}[1]{\widetilde {#1}}
\new{\bari}[1]{{#1}^{(i)}}
\new{\barun}[1]{{#1}^{(1)}}
\new{\barn}[1]{{#1}^{(n)}}
\new{\quot}[2]{\raise 2pt\hbox{#1}/\raise-4pt\hbox{#2}}
\new{\xb}{\bar x}
\new{\yb}{\bar y}
\new{\zb}{\bar z}

\new{\ga}{\gamma}

\newcommand{\s}{\sigma}

\newcommand{\pP}{{\mathbb P}}
\newcommand{\R}{{\mathbb R}}
\newcommand{\C}{{\mathbb C}}
\newcommand{\Q}{{\mathbb Q}}

\newcommand{\N}{{\mathbb N}}
\newcommand{\A}{{\mathbb A}}
\newcommand{\LND}{{\rm LND}}
\newcommand{\ML}{{\rm ML}}

\newcommand{\G}{{\Gamma}}

\newcommand{\p}{{\partial}}

\newtheorem{thm}{Theorem}[section]
\newtheorem{cor}[thm]{Corollary}
\newtheorem{lem}[thm]{Lemma}
\newtheorem{conj}[thm]{Conjecture}
\newtheorem{conv}[thm]{Convention}
\newtheorem{prop}[thm]{Proposition}
\newtheorem{nota}[thm]{Notation}
\newtheorem{defi}[thm]{Definition}

\newcommand{\bdefi}{\begin{defi}}
\newcommand{\edefi}{\end{defi}}
\newcommand{\bcor}{\begin{cor}}
\newcommand{\ecor}{\end{cor}}
\newcommand{\blem}{\begin{lem}}
\newcommand{\elem}{\end{lem}}
\newcommand{\bconv}{\begin{conv}}
\newcommand{\econv}{\end{conv}}
\newcommand{\bconj}{\begin{conj}}
\newcommand{\econj}{\end{conj}}
\newcommand{\bprop}{\begin{prop}}
\newcommand{\eprop}{\end{prop}}
\newcommand{\bthm}{\begin{thm}}
\newcommand{\ethm}{\end{thm}}
\newcommand{\bnota}{\begin{nota}}
\newcommand{\enota}{\end{nota}}
\newcommand{\be}{\begin{eqnarray}}
\newcommand{\ee}{\end{eqnarray}}
\newcommand{\no}{\noindent}

\title{Miyanishi's characterization of the affine 
3-space does not hold in higher dimensions}

\author{Shulim Kaliman}
\address{Department of Mathematics,
University of Miami,
Coral Gables, FL  33124, U.S.A.} 
\email{kaliman@@math.miami.edu}

\author{Mikhail Zaidenberg} 
\address{Universit\'e
Grenoble I, Institut Fourier, UMR 5582 CNRS-UJF, BP 74,
38402 St.\ Martin
d'H\`eres c\'edex, France} 
\email{zaidenbe@@ujf-grenoble.fr}

\thanks{Research of the  first author
partially supported by NSA grant  MDA904-00-1-0016.\\
This work was started during the second author's stay at the IHES; 
he would like  to thank IHES for the hospitality. \\
\mbox{\hspace{11pt}}{\it 1991 Mathematics Subject Classification}: 
13A02, 13F20, 14F15, 14J35, 14L30.\\
\mbox{\hspace{11pt}}{\it Key words}: affine space,
polynomial algebra, polynomial curve, locally nilpotent derivation, exotic structure} 

\begin{document}

\begin{abstract}  
We present an example which confirms the assertion of the title.
\end{abstract}

\maketitle 

\section*{Introduction}
Let $X$ be a smooth, contractible complex affine 3-fold. Recall

\bigskip \no {\bf Miyanishi's Theorem} 
\cite{Miy}\footnote{In the original formulation, 
instead of assuming $X$ to be topologically contractible, 
it is subjected to the following weaker conditions:
$e(X)=1$, the algebra $\C[X]$ of regular functions on $X$ is a UFD, and
all its invertible elements are constants. 
Likewise, in condition (ii) resp., (ii$'$) below the fibers themselves are replaced by their irreducible components. 
But actually, one can show that 
all the fibers of the function $f$ as in (i) 
are reduced and irreducible.}.
{\it $X\simeq \C^3$ if and only if the following two conditions hold:

\smallskip \no (i) there exists
a regular function $f\,:\,X\to\C$ and a Zariski open subset $U\subset \C$
such that $f^{-1}(U)\simeq_U U\times \C^2$} (in particular, the 
general fiber
$F_c:=f^*(c)\,\,\,(c\in U)$ of $f$ is isomorphic to the affine plane 
$\C^2$), {\it and 

\smallskip \no (ii) all the fibers $F_c\,\,\,(c\in\C)$
are UFD-s} (that is, for any $c\in\C$ the divisor 
$F_c$ is reduced and irreducible, and the algebra $A_c:=\C[F_c]$
of regular functions on the surface $F_c$ is a UFD). 

\medskip \no By   
\cite[Lemmas I, III]{Ka2} and \cite{KaZa2}, {\it the theorem holds  
if one only supposes that 

\smallskip \no ($i'$) the general fibers $F_c$ of $f$ 
are isomorphic to $\C^2$ and 

\smallskip \no ($ii'$)
each fiber $F_c\,\,\,(c\in\C)$ has at most isolated singularities.}

\smallskip \no The latter assumption ($ii'$) is essential, as shows 
the example of Russell's cubic 3-fold $X\subset \C^4$, $X=p^{-1}(0)$ 
where $p=x+x^2y+z^2+t^3$. 
In this example the fibers $F_c\,\,\,(c\in\C)$ of
the regular function $f=x\,|\,X\,:\,X\to\C$ 
are isomorphic to $\C^2$ except for the fiber $F_0$ which 
has non-isolated singularities (and therefore, it is not a UFD). 
And indeed, the Russell cubic $X$ is not isomorphic to $\C^3$ \cite{ML1},
that is, it is an {\it exotic $\C^3$} (i.e., a 
smooth affine variety diffeomorphic to $\R^6$ 
and non-isomorphic to $\C^3$;
see \cite{Za2}). 

More generally, the Main Theorem of 
\cite{Ka2, KaZa2} provides the following 
useful supplement to Miyanishi's theorem:
 
\medskip\no {\it A smooth, contractible affine 3-fold 
$X$ is an exotic $\C^3$
if there exists a regular function $f\,:\,X\to\C$ on $X$ with
general fibers isomorphic to $\C^2$, but not
all of its fibers being so.}\footnote{Notice that such a threefold $X$ 
admits a birational dominant morphism $\C^3\to X$.}

\medskip In this paper we prove the following

\bigskip

\noindent {\bf Theorem 1.} {\it  
The hypersurface $X$ in $\C^5$ given by the equation 
\be \label{th} p=u^mv+{(xz+1)^k-(yz+1)^l+z\over z}=0\,\ee
where $m\ge 2, \,\,k>l\ge 3,\,\,\,{\rm gcd}\,(k,l)=1$,
is an exotic $\C^4$.}

\medskip

\no {\bf Corollary.} {\it Miyanishi's Theorem does not hold in 
the dimension four.}

\medskip\no {\it Proof.}  The regular function 
$f=u\,|\,X\,:\,X\to\C$ on this hypersurface
provides a fibration  
with all the fibers $F_c\,\,(c\in\C)$ being smooth reduced
contractible affine 3-folds, all but the zero one $F_0$ 
being isomorphic to $\C^3$.
Moreover, the mapping 
$$(x,y,z,u)\longmapsto (x,y,z,u,\,v=u^{-m}q_{k,l}(x,y,z))$$
where
$$q_{k,l}:={(xz+1)^k-(yz+1)^l+z\over z}$$
gives an isomorphism $f^{-1}(U)=X\setminus F_0\simeq_U U\times \C^3$
where $U=\C\setminus\{0\}$. 
At the same time, the fiber $F_0\simeq S_{k,l}\times\C$ 
is an exotic $\C^3$ 
(see \cite{Za1, Za2}). Here $S_{k,l}=q_{k,l}^{-1}(0)\subset\C^3$
is the {\it tom Dieck-Petrie surface}; this  
is a smooth contractible affine surface non-homeomorphic to $\R^4$ 
\cite{tDP}. By
Fujita's theorem \cite[(1.18)-(1.20)]{Fu}, \cite[(3.2)]{Ka1}, 
any smooth, 
contractible affine variety is a UFD. Hence all the fibers 
of the function 
$f=u\,|\,X$ are smooth UFD-s. Thus the both conditions (i) 
and (ii) of the Miyanishi Theorem are fulfilled, whereas 
due to Theorem 1, $X\not\simeq\C^4$.
\qed

\medskip\no {\bf Remark.} Theorem 1 still holds 
for a triplet $(k,l,m)$ with
$l=2$ if ${\rm gcd}\,(m,2k)=1$. 
The proof of this fact is not difficult but we prefer the 
argument below since this enables us to demonstrate a nice
connection with the
Diophantine geometry over function fields (see Section 2). 
However, we do not know if the
statement remains true for (say) the triplet $(k,l,m)=(3,2,2)$. 

\medskip 
The proof of Theorem 1 is divided in two parts. 
The first one, concerning the topology of the variety $X$, 
is done in \cite{KaZa1}. The second one 
(which is done in Section 1 below) concerns exoticity of $X$; 
it mainly relies on the fact 
that there are only few regular actions on $X$ 
of the additive group $\C_+$ of the complex number field 
and moreover, there are only few polynomial curves in certain 
affine varieties related to $X$. 

\medskip
To conclude, recall the following

\medskip\no {\bf Problem.} {\it Let $X$ be a smooth, 
contractible 
complex affine
$n$-fold where $n\ge 4$, and let $f\,:\,X\to\C$ be 
a regular function on $X$. 
Suppose that $f^*(c)\simeq \C^{n-1}$ for every  
$c\in\C$. Is it true that $X\simeq\C^n$, and that 
this isomorphism sends 
$f$ into a variable of the polynomial algebra $\C^{[n]}$?}

\medskip\no The results of \cite{Miy, Ka2, KaZa2, Sa} 
cited above provide 
a positive answer for\footnote{It is worthwile 
noting that, without the assumption that $X$ is affine,
the answer is negative even for 
$n=4$.
Indeed, consider the smooth non-affine 4-fold 
$X={\tilde X}\setminus Z$ where ${\tilde X}$ 
is the hypersurface $uv=xy + z^2-1$ in $\C^5$ and $Z$ is 
the plane $u=x=z-1=0$ in ${\tilde X}$. Then every fiber
of the morphism $(x,u)\,:\,X\to\C^2$ is isomorphic to 
$\C^2$, and every fiber of the regular function $u\,|\,X$
is isomorphic to $\C^3$.}
$n=3$. 

We are grateful to the referee for useful 
remarks and suggestions
which served us to improve the exposition. 

\section{Proof of Theorem 1}

In \cite[Proposition 4.4, Example 6.2]{KaZa1}
it is shown that the smooth affine 4-fold $X$ as 
in Theorem 1 is contractible and moreover, 
diffeomorphic to $\R^6$.
Thus, to prove the theorem
it is enough to verify that $X\not\simeq \C^4$. 
The proof of the latter assertion is based 
on the computation of the
{\it Makar-Limanov
invariant} $\ML(X)$. Recall that $\ML(X)$ 
denotes the algebra 
of regular functions
on the variety $X$ invariant under any regular 
$\C_+$-action on $X$ (or in
other words, of regular functions on $X$
that are vanished by any locally nilpotent
derivation of the algebra $\C[X]$; see e.g., 
\cite{KaML1, Za2}, or also \cite{De}).

In fact, we prove
the following

\medskip\no {\bf Proposition 1.} 
$\ML(X)\supset \C[\ol  u]$ where $\ol u = u\,|\,X$. 
{\it Hence 
$\ML(X)\not\simeq \ML(\C^4)=\C$ and therefore, 
$X\not\simeq\C^4$.}

\medskip\no {\bf Remark.} If $m=1$ then $\ML(X)=\C$ 
(and moreover, the  group of biregular automorphisms 
of $X$ generated by the regular $\C_+$-actions 
on $X$ acts infinitely transitively 
\cite[Theorem 5.1]{KaZa1}). The question arises: 
{\it is it still true that $X$ is an exotic 
$\C^4$ when $m=1$, at least for some values of $k$ and $l$? }

\bigskip

\no {\bf Notation.} Throughout
the proof, we fix a weight degree function $d$ on the 
polynomial algebra $\C^{[5]}=\C[x,y,z,u,v]$
given by 
\be\label{we} d_x=l,\quad d_y=k,\quad d_z=0,\quad 
d_u=-n\sqrt{2},\quad d_v=mn\sqrt{2}+kl\,\ee 
where $n\in\N$.
This degree function $d$
satisfies the following conditions:
$$kd_x+(k-1)d_z=ld_y+(l-1)d_z=md_u+d_v = kl$$
$$>\max_{i=1,\dots,k-1,\,\,j=1,\dots,l-1} \{0,\,\,id_x+(i-1)d_z=il,\,\,jd_y+(j-1)d_z=jk\}\,.$$
It follows that 
$$\h p:=u^mv+x^kz^{k-1}-y^lz^{l-1}= 
u^mv+z^{l-1}(x^kz^{k-l}-y^l)$$
is the principal $d$-homogeneous part of 
the polynomial $p$ from
(\ref{th}); indeed, 
\be\label{po} 
p=u^mv+q_{k,l}(x,y,z)=u^mv+\sum_{i=1}^k 
{k \choose i}x^iz^{i-1}-
\sum_{j=1}^l {l\choose j}y^jz^{j-1}+1\,.\ee
By  $d_A$ we denote the induced degree function 
on the
algebra $A:=\C[X]=\C^{[5]}/(p)$. Let $\h A$
be the associate graded algebra, and $d_{\h A}$ 
be the induced 
degree function on $\h A$.
Since the polynomial $\h p$ is irreducible, by 
Proposition 4.1 in \cite{KaML3} 
(see also \cite[Lemma 7.1]{Za2}),
the affine variety $\h X:=$spec$\,\h A$ 
coincides with the hypersurface
in $\C^5$ given by the equation $\h p=0$.
We denote by
$\h x\dt \h v$ the images in $\h A$
of the coordinate functions $x\dt v$,  
respectively,
whereas their restrictions to $X$ are 
denoted as $\ol x\dt \ol v$. Thus
in the algebra $\h A$ the following 
relation holds:
\be\label{re}\h u^m\h v=\h z^{l-1}
(\h y^l-\h x^k\h z^{k-l})\,.\ee

For an integral domain $B$ of finite type, 
let  $\LND(B)$ be the set of all its 
locally nilpotent derivations. 
Fix arbitrary $\p\in\LND(B)\setminus \{0\}$.
Recall the following well known facts 
which we frequently use below
(see e.g., \cite{ML1, KaML1, Za2}).

\bigskip

\noindent {\bf Lemma 0.}  {\it
(a)
The invariant subalgebra ${\rm ker}\,\p \subset B$
is factorially closed, that is, 
$ab\in {\rm ker}\,\h \p \setminus \{0\} 
\Longrightarrow a,\,b\in {\rm ker}\,\h \p$. 
Moreover,\footnote{The latter statement 
is a lemma due to Makar-Limanov
(see e.g., \cite[Ex.(7.12.d)]{Za2}), 
which also follows from the Davenport Lemma 
in Section 2 below.} 
$$a^k+b^l\in {\rm ker}\,\h \p \setminus \{0\} \quad 
\mbox{and} \quad k,\,l\ge 2 
\Longrightarrow a,\,b\in {\rm ker}\,\h \p\,.$$

\smallskip\noindent (b) Let $a\in B$ be an element
of $\p$-degree one, i.e., 
$\p a \in {\rm ker}\,\h \p \setminus \{0\}$. 
Then
any element $b\in B$ can be presented in the form
\be \label{da} b=c^{-1}\sum_{i=0}^N c_ia^i\,\ee 
where $c,\,c_0,\dots,c_N\in {\rm ker} \,\p$. 

\smallskip\noindent (c) 
The invariant subfield ${\rm Frac \,ker}\,\p \subset {\rm Frac }\,B$ 
is algebraically closed in the fraction field 
${\rm Frac}\,B$, 
and
${\rm tr.deg \,[Frac}\,B : {\rm Frac\, ker} \,\p]=1$. }

\bigskip

Fix a locally nilpotent derivation 
$\p \in \LND(A)\setminus \{0\}$, 
and
let $\h \p\in \LND(\h A)$ be 
the homogeneous locally nilpotent derivation 
of the graded algebra $\h A$ associated with $\p$
(that is, the principal part of $\p$);  
notice that $\h \p \neq 0$ once $\p \neq 0$
(see \cite{ML1} or also \cite{KaML1, Za2}).

\bigskip

\noindent {\bf Lemma 1.} ker$\,\h \p \not\subset 
\C[\h x,\h y,\h z]$.

\medskip

\noindent {\it Proof.} Assume the contrary. 
Since tr.deg$\,[\h A :$ ker$\,\h \p ]=1$ 
there exist three algebraically independent 
elements, say,
$a,\,b,\,c\in$ ker$\,\h \p $. Regarding 
the elements $a,\,b,\,c$ 
as polynomials
in $x,y,z$, consider the morphism 
$\s=(x,a,b,c)\,:\,\C^3\to\C^4$.
The Zariski closure of the image of $\s$ 
being a proper algebraic subvariety of
$\C^4$, there is a non-trivial relation 
$g(x,a,b,c)=0$ where
$g\in\C^{[4]}\setminus \{0\}$. Hence 
we have a non-trivial relation
\be\label{g} \sum_{i=0}^N g_i(a,b,c)\h x^i=0\ee 
in the algebra $\h A$ where $g_i(a,b,c)\in $ 
ker$\,\h \p$. Here
$N>0$ (indeed, otherwise the elements $a,\,b,\,c$ 
would be algebraically 
dependent). It follows from Lemma 0(c) that 
$\h x\in \,$ker$\,\h \p $. 
Similarly, we have $\h y,\,\h z\in $ker$\,\h \p $. 
In virtue of the relation
(\ref{re}) above, also $\h u^m\h v\in$ ker$\,\h \p $. 
By Lemma 0(a), 
it follows that 
$\h u,\,\h v \in {\rm ker}\,\h \p $. Therefore,
ker$\,\h \p =\C[\h x,\h y,\h z,\h u,\h v]=\h A$, 
and so 
$\h \p =0$, a contradiction. 
This proves the lemma. \qed

\medskip

\noindent {\bf Lemma 2.} {\it The following alternative holds:
either 
$\h u\in {\rm ker}\,\h \p$ or $\h v\in {\rm ker}\,\h \p$.}

\medskip

\noindent {\it Proof.} Due to (\ref{re}), 
any element $\h a\in\h A=\C[\h X]$ can be extended to a unique polynomial
$\h f\in\C^{[5]}$ of the form
\be\label{expr} 
\h f=\sum_{i\ge 0}a_iu^i+\sum_{i=0}^{m-1}\sum_{j>0}b_{ij}u^iv^j\,\ee
where $a_i,\,b_{ij}\in\C[x,y,z]$. It is known \cite{KaML3} 
(see also \cite[Ex.(7.8)]{Za2}) that 
$$d_{\h A}(\h a)\stackrel{def}{=}\min\{d(f)\,|\,f\in\C^{[5]},
\,\,f\,|\,\h X=\h a\}=\min\,\{d(f)\,|\,f=
\h f+\h pg,\,\,g\in\C^{[5]}\}=d(\h f)\,.$$ 
Furthermore, if $\h a\in\h A$ is a
$d_{\h A}$-homogeneous element, then the polynomial $\h f$ is
$d$-homogeneous, too.

Since the derivation $\h \p $ is homogeneous (i.e., graded) its kernel
ker$\,\h \p $ is a graded subalgebra of the graded algebra $\h A$, and so 
it is generated by homogeneous elements. Let $\h a\in$ ker$\,\h \p $ 
be a non-zero homogeneous element, and let 
$\h f\in \C^{[5]}$ be its $d$-homogeneous extension as in (\ref{expr}) above. 

The degree function $d\,:\,\C^{[5]}\setminus\{0\}\to \Q[\sqrt{2}]$ 
can be represented as $d=d'+\sqrt{2}d''$ where 
$d',\,d''\,:\,\C^{[5]}\setminus\{0\}\to \Q$.
By (\ref{we}), we have 
$d\,|\,\C[ x,y,z]=d'\,|\,\C[x,y,z]$ and
$$d''(a_iu^i)=-in,\quad d''(b_{ij}u^iv^j)=
(jm-i)n,\,\,i=0,\dots,m-1\,$$
assuming that $a_i,\,b_{ij}\neq 0$.
All these degrees are pairwise distinct. Hence the polynomial  $\h f\neq 0$ 
being $d$-homogeneous, the 
expression (\ref{expr}) for $\h f$ consists of a single term. 
Therefore, the following alternative holds: 

\smallskip

\no (i)  either $\h a=a_0\in\C[\h x,\h y,\h z]$, or

\smallskip

\no (ii) $\h a=a_i\h u^i$ for some $i>0$ and for some 
$a_i\in\C[\h x,\h y,\h z]\setminus\{0\}$, or

\smallskip

\no (iii) $\h a=b_{ij}\h u^i\h v^j$ for some $j>0$ and for some 
$b_{ij}\in\C[\h x,\h y,\h z]\setminus\{0\}$.

\smallskip

\no Since the subalgebra ker$\,\h \p $ is factorially closed, 
we have $\h u\in$ ker$\,\h \p $ in the case (ii) 
and $\h v\in$ ker$\,\h \p $ in the case (iii). 
By Lemma 1, (i) cannot happen for all the homogeneous
elements $\h a\in$ ker$\,\h \p$. Thus, the assertion follows. \qed

\medskip

\noindent {\bf Lemma 3.} $\h v\not\in {\rm ker}\,\h \p $.

\medskip

\noindent {\it Proof.} Assume on the
contrary that $\h v\in$ ker$\,\h \p $.
Then for a general $c\in\C$, the locally nilpotent derivation $\h \p $
can be specialized to a locally nilpotent derivation 
$\h \p _c\in\LND(\h X_c)\setminus \{0\}$ 
where for $c\in\C\setminus\{0\}$ we denote
$$\h X_c=\h X\cap \{v=c\}\simeq 
\h X_1=\{u^m+z^{l-1}(x^kz^{k-l}-y^l)=0\}\subset\C^4\,.$$ 
We keep the same notation $\h \p $ for $\h \p _1$, and we still denote by 
$\h \varphi$ the associated $\C_+$-action $\h \varphi\,|\,\h X_1$ on the threefold $\h X_1$. 

Note that the threefold $\h X_1$ has divisorial singularities.
Indeed, since by our assumption, $m\ge 2$ and $k>l\ge 3$, it is singular along
the divisor $D_{\h z}$ of the regular function $\h z\in \C[\h X_1]$:
$D_{\h z} \subset {\rm sing}\,\h X_1$. It follows that the divisor $D_{\h z}$
is invariant under the $\C_+$-action $\h \varphi$ on $\h X_1$.
Hence a general $\h \varphi$-orbit $\cal O$ does not meet the divisor
$D_{\h z}$, and so the restriction $\h z\,|\,{\cal O}$ does not
vanish. 
Therefore, the regular function $\h z$ is constant along general $\h \varphi$-orbits, that is,
$\h z$ is a $\h \varphi$-invariant, or equivalently,
$\h z\in{\rm ker}\,\h \p$. 

Thus, we are in the position to repeat the specialization descent. Namely,
the $\C_+$-action $\h \varphi$ can be further specialized to the 
general $\h \varphi$-invariant surface 
$$S_c:=\{\h z=c\}\simeq S_1=\{u^m+x^k-y^l=0\}\subset \C^3\,$$
providing a non-trivial $\C_+$-action on $S_c$ and thereby also on $S_1$. 
Now the desired conclusion follows from the next lemma.
\qed

\medskip

\no {\bf Lemma 4.} {\it The Pham-Brieskorn surface 
$$S=S_{k,l,m}=\{x^k+y^l+z^m=0\}\subset\C^3$$
where $k,l,m\ge 2$  
admits a non-trivial regular $\C_+$-action if and only if this is a dihedral surface
$S_{2,2,m}$. In the latter case ML$\,(S)=\C$.}

\medskip

\no {\it Proof.} Assume that $\h \varphi$ is a non-trivial regular 
$\C_+$-action on $S$. Let ${\cal O}\subset S$ 
be a general $\h \varphi$-orbit. Since ${\cal O}\simeq \C$,
it can be parameterized by a triple of polynomials
$(x(t),\,y(t),\,u(t))\in (\C[t])^3$ 
satisfying the relation
$$x^k(t)+y^l(t)+z^m(t)=0\,.$$ 
Assume first that $1/k+1/l+1/m\le 1$. 
Then by the Halphen Lemma (a)
in the next section (these
polynomials cannot be relatively prime in pairs, and so) 
a general orbit $\cal O$ 
meets one of the axes $x=y=0$ or $x=z=0$ or $y=z=0$, hence it
must pass through the origin, a contradiction. 

In the remaining cases
$S$ is one of the Platonic surfaces 
$S_{2,2,m},\,\, S_{2,3,3},\,\,S_{2,3,4}$ or $S_{2,3,5}$. 
Anyhow, to exclude the last three cases
we will assume in the sequal more generally that 
${\rm gcd}\,(m,kl)=1$, and that on the contrary,
$\LND(S)\neq\{0\}$, that is, that 
the surface $S$ admits a non-trivial regular $\C_+$-action. 

Let $\p_0\in\LND(S)$, $\,\p_0\neq 0$. Fix a weight degree
function  $d'$ on the polynomial algebra $\C^{[3]}$ given by
$d'_x=1/k,\,\,d'_y=1/l,\,\,d'_z=1/m$. Since the polynomial
$x^k+y^l+z^m$ is $d'$-homogeneous, the algebra $B=\C[S]$ is graded. 
The graded locally nilpotent derivation $\h \p_0$ of the algebra $B$
associated with $\p_0$ is also non-zero. In virtue of the relation 
$\h z^m=-(\h x^k+\h y^l)$, any element $\h b\in B$ extends to a unique
polynomial $\h f\in\C^{[3]}$ with $\deg_z\,\h f<m$. If the element
$\h b$ is $d'_B$-homogeneous then also the
polynomial $\h f$ is $d'$-homogeneous, and 
the following statement holds.

\medskip\no {\bf Claim.} $\h f=cx^{\alpha}y^{\beta}z^{\gamma}\prod_i
(x^{k'}-c_iy^{l'})$ {\it where $k'=k /{\rm gcd}\,(k,l),\,\,l'=l /{\rm gcd}\,(k,l),\,\,c,\,c_i\in\C^*$, and $\gamma<m$.} 

\medskip\no  {\it Proof of the claim.}  Letting
$d'(x^iy^jz^s)=d'(x^{i'}y^{j'}z^{s'})$ where $0\le s\le s'\le m-1$ we
will have \be\label{7} {s'-s\over m}= {i-i'\over k}+{j-j'\over l}\,.\ee
Since by our assumption ${\rm gcd}\,(m,kl)=1$, it follows from  (\ref{7}) that $m\,|\,(s'-s)\Longrightarrow s=s'$, and hence
${i-i' \over k}= {j'-j \over l}\Longrightarrow {i-i' \over k'}= {j'-j \over l'}$. Now the claim follows.  \qed

\medskip

The graded subalgebra ${\rm ker}\,\h \p_0$ of the algebra $B$ being
generated by homogeneous elements, there exists a non-zero homogeneous
element
$b\in {\rm ker}\,\h \p_0$. Since the subalgebra ${\rm ker}\,\h \p_0$
is factorially
closed, in virtue of the above claim, the following
alternative holds:

$${\rm (i)}\quad \h x\in {\rm ker}\,\h \p_0\quad {\rm or}\quad 
{\rm (ii)}\quad  \h y\in {\rm ker}\,\h \p_0\quad {\rm or}\quad
{\rm (iii)}\quad \h z\in {\rm ker}\,\h \p_0\quad {\rm or}$$ 
$${\rm (iv)}\quad  \h x^{k'}-c_i\h y^{l'}\in {\rm ker}\,\h \p_0\qquad {\rm 
 for \,\,\,some}\qquad c_i\neq 0\,.$$ 

\no Since $\h x^k+\h y^l+\h z^m=0$, in the case (i) we have 
$\h y^l+\h z^m\in {\rm ker}\,\h \p_0$. As $l,\,m\ge 2$ 
by Lemma 0(a) this implies
$\h y,\,\h z\in  {\rm ker}\,\h \p_0$.  
Henceforth, $\h \p_0=0$, a contradiction. Similarly, the cases (ii)
and (iii) lead to  a contradiction. 

If $\min\,\{k',l'\}\ge 2$ then by the same arguments as above, (iv) implies 
that $\h x,\,\h y\in  {\rm ker}\,\h \p_0$, and then also 
$\h z\in  {\rm ker}\,\h \p_0$, which again gives a contradiction. Thus it must be $\min\,\{k',l'\}= 1$; let $l'=1$. Then $k=lk'$. The regular function 
$\h x^{k'}-c_i\h y \in\C[S]$ being  invariant under the associated regular $\C_+$-action $\varphi_{\h \p_0}$ on the surface $S$, its general level sets
contain general $\varphi_{\h \p_0}$-orbits. Being irreducible, these curves should be isomorphic to $\C$. On the other hand, they are isomorphic to the affine plane curves with the equations
$$x^{lk'}+\left({x^{k'}-c'\over c_i}\right)^l+z^m=0$$
where $c'\in\C$ is generic. It is easily seen that such a curve cannot be isomorphic to  $\C$ unless $k=l=2$ and $c_i^2=-1$, in which case $S$ is a dihedral surface (hint: notice that 
an irreducible affine curve is isomorphic to $\C$ if and only if it admits a regular $\C_+$-action, and then proceed in the same fashion as above). 

To prove the last statement of the lemma, notice that
there is an isomorphism
$S_{2,2,m}\simeq T_m:=\{uv-w^m=0\}$, 
and hence ML$(S_{2,2,m})\simeq$ ML$(T_m)=\C$. 
The latter equality is well known; 
see e.g. \cite{DanGi, Be, ML2, ML3, KaZa1}. 
Indeed, the subgroup $\langle\alpha,\,\beta\rangle$ 
of the automorphism 
group Aut$\,T_m$ generated by the following $\C_+$-actions on $T_m$ 
(restricted from $\C^3$):
$$\alpha\,:\,(t,\,(u,v,w))\longmapsto 
\left(u,\,v+{(w+tu)^m-w^m\over u},\,w+tu\right)\,,$$ 
$$\beta\,:\,(t,\,(u,v,w))\longmapsto 
\left(u+{(w+tv)^m-w^m\over v},\,v,\,w+tv\right)\,.$$ 
has a dense orbit; 
therefore, ML$\,(T_m)=\C$. 
This concludes the proof. \qed

\medskip

From Lemmas 1-3 we obtain such a corollary.

\medskip

\noindent {\bf Corollary.} 
$\h u \in {\rm ker}\,\h \p  $.

\medskip

\no {\bf Lemma 5.} 
${\h \p } \h v \not\in {\rm ker}\,\h \p \,.$

\medskip

\no {\it Proof.} 
Assume the contrary. 
Then by (\ref{re}) and the above corollary,
we have:
$$\h \p  (\h u^m \h v)= \h u^m\h \p  
\h v \in {\rm ker}\,\h \p 
\Longrightarrow \h \p  [\h z^{l-1} 
(\h y^l - \h x^k \h z^{k-l})]=\h \p  
(g_1g_2) \in
{\rm ker}\,\h \p \,$$
where $g_1:=\h z^{l-1}$ and 
$g_2:=\h y^l - \h x^k \h z^{k-l}\in\C[\h X]$.
Hence the restriction of 
the product $g_1g_2$ onto 
a general orbit $\cal O$
serves as a coordinate function 
of the curve ${\cal O}\simeq \C$
(for instance, this follows from Lemma 0(b)). 
In other words, 
$\deg_t  [(g_1g_2)\,|\,{\cal O}]=1$ 
where 
$t$ is a coordinate in ${\cal O}\simeq \C$
(notice that $\deg_t  (f\,|\,{\cal O})=\deg_{\h \p} f$ 
where the latter degree is defined below). 
This provides the following alternative: 

\smallskip \no $\bullet$ either 
$\deg_t  (g_1\,|\,{\cal O})=0$ 
and $\deg_t  (g_2\,|\,{\cal O})=1$, or

\smallskip \no $\bullet$ 
$\deg_t  (g_1\,|\,{\cal O})=1$ 
and $\deg_t  (g_2\,|\,{\cal O})=0$.

\smallskip\no 
Consider each of these two possibilities.

\smallskip\no Assuming first that 
$z\,|\,{\cal O}={\rm const}\in\C\setminus\{0\}$ 
(i.e., $\h z \in {\rm ker}\,\h \p $) and 
$\deg_t  [(y^l-x^kz^{k-l})\,|\,{\cal O}]=1$, 
we would have that $\deg_t  (y^l(t)-cx^k(t))=1$ 
for two polynomials $x(t),\,y(t)\in\C[t]$
and for a general constant $c\neq 0$. 
We may also suppose that ${\rm
  gcd}\,(x(t),\,y(t))=1$
(i.e., that the orbit $\cal O$ 
does not meet the codimension two
subvariety $D_{\h x} \cap D_{\h y}$ of $\h X$). 
Then by the Davenport Lemma
in the next section, for a certain 
$m'\in\N $ the inequalities
$1>m'(kl-k-l)\ge 1$ must hold, 
which is impossible. 

In the second case we would have:
$\deg_t  (g_1\,|\,{\cal O})=
(l-1)\deg_t (z\,|\,{\cal O}) = 1$
which is also impossible since 
by the assumption of Theorem 1, 
$l-1\ge 2$.
This proves the lemma. \qed

\medskip Recall the notion of 
the degree function associated 
with a locally nilpotent derivation 
$\p\in\LND(A)$:
$$\deg_{\p}\,a\stackrel{def}{=}
\max\,\{n\in\N\cup\{0\}\,|\,\p^n a
\neq 0\}\quad{\rm if}\quad a\in 
A\setminus \{0\};\,\,\,\,\deg_{\p}\,0=
-\infty\,.$$
For the associated locally nilpotent derivation
$\h \p\in\LND(\h A)$ we have the inequality
$$\deg_{\p} a \ge \deg_{\h \p} \h a 
\qquad\forall a\in A\,,$$
where $\h a\in\h A$ denotes 
the principal $d$-homogeneous part of $a$.

Lemma 5 provides the following

\medskip\no {\bf Corollary.} 
$ \p v \not\in {\rm ker}\,\p \,;$ 
{\it moreover, $\deg_{\p} v \ge \deg_{\h \p} \h v\ge 2$.}

\medskip\no {\bf Lemma 6.} 
{\it Let $a\in A$ be an element such that
$\deg_{\p} a \le 1$. Then $a$ 
can be extended to a polynomial 
$f\in\C^{[5]}$
which does not depend on $v$.}

\medskip\no {\it Proof.} Since 
$p\,|\,X=0$, in virtue of (\ref{po}) the
restriction $\ol u^m\ol v$ of the polynomial 
$u^mv$ to $X$ can be expressed
as a polynomial in $\ol x,\,\ol y$ and $\ol z$. 
Hence the element $a\in A$
can be extended (in a unique
way) to a polynomial $f\in\C^{[5]}$ 
written in the form (\ref{expr}).
Let us show that this polynomial $f$ 
does not depend on $v$.
Assume the contrary. 
Letting the constant $n$
in the definition (\ref{we}) 
of the weight degree
function $d$ be large enough, 
we can achieve
(by the same arguments as 
in the proof of Lemma 2) that
the principal $d$-homogeneous part 
$\h f$ of $f$ is as in (iii) 
of this same proof. In particular, 
$v$ is a factor of the polynomial
$\h f$. Hence $\h v$ is a factor 
of $\h a =\h f\,|\,\h X$. By Lemma 5, 
we have the inequalities
$$\deg_{\p} a \ge \deg_{\h \p} \h a 
\ge \deg_{\h \p} \h v \ge 2$$
which contradicts our choice of the element $a$. 
The lemma is proven.
\qed

\medskip\no {\it Proof of Proposition 1} 
(cf. \cite{KaML1, ML3}). 
We have to show that 
$u\in {\rm ker} \,\p$ for any $\p \in \LND(A)$. 
Fix an element $a\in A$ 
of $\p$-degree one. 
Letting in (\ref{da}) $b=v$, 
from (\ref{po}) and (\ref{da}) we obtain:
\be\label{8}\ol v=c^{-1}\sum_{i=0}^N c_ia^i=
-{\ol u}^{-m}q_{k,l}(\ol x,\ol y,\ol z)\quad
\Longrightarrow \quad-\ol u^m \sum_{i=0}^N 
c_ia^i= cq_{k,l}(\ol x,\ol y,\ol
z)\,.\ee
By Lemma 6, the element $c\in A$ 
resp., 
$\sum_{i=0}^N c_ia^i\in A$ can be
extended to a polynomial, say, 
$\eta\in\C[x,y,z,u]$ resp., 
$\zeta\in\C[x,y,z,u]$. By
(\ref{8}), there exists a polynomial 
$g\in\C^{[5]}$ such that 
\be\label{9} u^m\zeta(x,y,z,u)-
q_{k,l}(x,y,z)\eta(x,y,z,u)=pg\,.\ee
The left hand side of (\ref{9}) 
does not depend on $v$ but the
polynomial $p$ does, hence we must 
have $g=0$. Since 
${\rm gcd}\,(u,\,q_{k,l})=1$
it follows from (\ref{9}) that $u$ 
divides $\eta$ in the algebra 
$\C^{[4]}$ and
so, $\ol u$ divides $c$ in the algebra $A$,
that is, $c=\ol u b$
where $b\in A$. Since $c\in {\rm ker}\, \p$ 
and $ {\rm ker}\, \p $ is
factorially closed, also 
$\ol u\in {\rm ker}\, \p$,
as stated. This completes the proof.
\qed

\section{$\A_1$-poor varieties: 
the lemmas of Mason, Davenport  
and Halphen}

In course of the proof of Proposition 1 
we have used the lemmas of 
Davenport and Halphen; 
for the sake of completeness, 
we provide them below with simple proofs
based on the following well known

\medskip\no {\bf Mason's abc-Lemma \cite{Mas}.} 
{\it Let $a,\,b,\,c\in\C[t]$ be three 
polynomials, not all three constant. 
For a polynomial $p\in\C[t]$, 
denote by $d_0(p)$ the number of its
distinct roots 
(without counting multiplicities). 
Assume that $a+b+c=0$ and
${\rm gcd}\,(a,b)=1$.
Then we have
\be\label{ma} \max\{\deg\, a,\,\deg\, b,\,\deg\, c\}
\le d_0(abc)-1\,.\ee}
 
\smallskip\no See \cite{La, Mas, Si} 
for an elementary proof. 
We would like to sketch 

\medskip \no {\it An alternative proof.} 
Let $f\,:\,\G_1\to\G_2$ be a  
proper, surjective morphism 
of smooth quasiprojective curves. 
Then the following inequality for 
Euler characteristics holds: 
\be\label{RH} e(\G_1)\le (\deg\,f)\,e(\G_2)\,.\ee
This inequality follows from the obvious relations
$${\rm card}\,({\rm CrPt(f)})
\le (\deg\,f)\,{\rm card}\,({\rm CrVa}(f))$$
and 
$$e(\G_1\setminus {\rm CrPt}(f))=
(\deg\,f)\,e(\G_2\setminus {\rm CrVa}(f))$$
where ${\rm CrPt}(f)$ resp., 
${\rm CrVa}(f)$ denotes the set of critical points
resp., critical values of $f$. 

Take $\G_1=R\setminus S$ where $R$ is 
a smooth projective curve of genus $g$ 
and $S$ is a finite subset of $R$, and 
let $\G_2\simeq\C\setminus\{0,1\}$ be 
realized as $\G_2=\{u+v=1,\,\,u
\neq 0,\,\,v\neq 0\}\subset\C^2$. 
Then for a pair $f=(u,v)$ of non-constant 
rational functions on $R$ with zeros and 
poles only on $S$
such that $u+v=1$,
from (\ref{RH}) we obtain the inequality 
(see \cite{Mas})
\be\label{mi} \deg\,u= \deg\,v\le-e(R\setminus S)
=2g-2+{\rm card}\,S\,.\ee
Letting 
$$R=\pP^1,\qquad S=\{\infty\}
\cup a^{-1}(0)\cup b^{-1}(0)\cup c^{-1}(0),\qquad u
=-a/c,\quad v=-b/c$$
(so that the condition  $a+b+c=0$ of the lemma 
becomes $u+v=1$), from (\ref{mi})
we get (\ref{ma}).\qed

\medskip
As an immediate corollary, we obtain

\medskip\no {\bf Davenport's Lemma 
\cite[Thm.2]{KlNe, Dav}\footnote{See also 
\cite{DvZa} and the literature therein 
for closely related results.}.} 
{\it Let three
polynomials $x,\,y,\,z\in\C[t]$  
satisfy the relation $z=x^k-y^l$ where
$k$ and $l$ are relatively prime\footnote{One 
can find in \cite{Dav} a general formulation 
with arbitrary $k$ and $l$.}, $z\neq 0$, 
${\rm gcd}\,(x,y)=1$  and 
$\deg z < \max\{\deg x^k,\,\deg y^l\}$.
Denote $n=\deg z,\,\,lm=\deg x,\,\,km=\deg y$. 
Then we have
$$n> m(kl-k-l)\,.$$}
\no {\it Proof} (cf. \cite{Pr}). 
By Mason's abc-Lemma, we have the inequality
$$\max\{k\deg x,\,l\deg y\} 
\le \deg x +\deg y +\deg z - 1\,.$$
Hence $$klm\le km+ lm+n-1\,,$$
and the lemma follows.\qed

\medskip\no {\bf Remark.} 
It is known \cite{St, Zn, Ore1} that 
(whatever $k,\,l$ and $m$ with 
${\rm gcd}\,(k,l)=1$ are)
the bound in Davenport's Lemma 
is the best possible one. 
See also \cite{Si} on exactness 
in Mason's abc-Lemma. 

\bigskip A contemporary exposition 
of Halphen's results \cite{Ha} 
is given in \cite{BaDw}. Actually, 
the original Halphen's Lemma has 
a broader meaning in the context of 
our subject. To formulate it in an 
appropriate way, we introduce the 
following notions\footnote{Cf. the 
notion of abc-variety in \cite{Bu}. 
Presumably (over the field $\C$) 
these are the affine varieties $X$ 
which do not admit non-constant morphisms 
$\C^*\to X$ \cite[p.231]{Bu}.}.

\medskip\no {\bf Definition.} Let 
$X$ be an algebraic variety.
We say that $X$ is {\it $\A_1$-poor} 
if 
there exists a subvariety $Y$ of $X$ 
of codimension at least two such that
every curve (i.e., a non-constant morphism) 
$f\,:\,\C\to X$  
meets $Y\,\,$: $\,\,\,f(\C)\cap Y 
\neq\emptyset$. 

\smallskip
In contrast, we say that $X$ is 
{\it $\A_1$-rich} if, 
for any two disjoint closed subvarieties 
$Y,\,Z\subset X$ with codim$_X Y\ge 2$ 
and dim$\,Z = 0$, there exists a 
polynomial curve $\C\to X$ omitting $Y$
and passing through every point of $Z$.

\medskip\no {\bf Remarks. 1.} 
Evidently, an $\A_1$-poor variety $X$ 
does not admit
non-trivial regular $\C_+$-actions. 
Or equivalently, 
LND$\,(X)=\{0\}$ $\Longleftrightarrow 
\ML(X)=\C[X]$. 
Moreover, the latter equality holds 
assuming that the algebra
$A=\C[X]$ is endowed with a degree 
function such that
for the associated graded algebra $\h A$, 
the variety $\h X={\rm spec}\,\h A$
is $\A_1$-poor. This justifies our 
interest in $\A_1$-poor varieties.

\smallskip\no {\bf 2.}  
The affine space $\C^n\,\,(n\ge 2)$ 
is $\A_1$-rich. 
Indeed, given two disjoint closed 
subvarieties 
$Y,\,Z\subset \C^n$ with 
codim$_{\C^n} Y\ge 2$ 
and dim$\,Z = 0$, 
by a theorem due to Gromov and 
Winkelmann \cite{Grm, Wi}, 
one can find an automorphism 
$\alpha\in {\rm Aut}\,\C^n$ 
such that 
$\alpha(Y)=Y$ and the image 
$\alpha(Z)$ is contained in an affine line
$L\subset \C^n\setminus Y$. 
Then the polynomial curve 
$\C\simeq L\stackrel{\alpha^{-1}}{\longrightarrow}
\alpha^{-1}(L)\subset \C^n$ omits $Y$
and passes through every point of $Z$, 
as required. 

\smallskip\no {\bf 3.} If a variety $X$ 
admits a finite morphism $X'\to X$ from 
an $\A_1$-rich affine variety $X'$
(for instance, from $X'=\C^n,\,\,n\ge 2$), 
then clearly $X$ is also $\A_1$-rich.
Notice also that the family of polynomial 
curves in an $\A_1$-rich affine variety is 
unbounded (that is, their degrees are not bounded 
in common). 

\medskip The following two lemmas provide 
examples of $\A_1$-poor resp., 
$\A_1$-rich surfaces in $\C^3$.

\medskip\no {\bf Halphen's Lemma 
\cite{Ha, Ev, BaDw}.} 
{\it Consider the Pham-Brieskorn surfaces
$$S_{k,l,m}=\{x^k+y^l+z^m=0\}\subset \C^3\,$$
where $k,l,m\ge 2$. Then the following 
statements hold.

\smallskip\no (a) The surface $S_{k,l,m}$
is $\A_1$-poor if and only if $1/k+1/l+1/m\le 1$. 
Actually, under the latter condition
any polynomial curve $f\,:\,\C\to S_{k,l,m}$ 
passes through 
the singular point
${\overline 0}\in S_{k,l,m}\subset \C^3$.

\smallskip\no (b) In contrast, every 
Platonic surface $S_{k,l,m}$
where $1/k+1/l+1/m> 1$ is $\A_1$-rich.}

\medskip\no {\it Proof.} (a) Suppose 
first that $1/k+1/l+1/m\le 1$. Let us 
show that no triple of non-constant 
relatively prime polynomials $(x(t),\,y(t),\,z(t))$ 
satisfies the relation
$x^k+y^l+z^m=0$. Assuming the contrary, 
by Mason's abc-Lemma, we have:
$$\max\{k\deg x,\,l\deg y,\,m\deg z\} 
\le \deg x +\deg y+\deg z  - 1\,.$$
Thus,
\begin{eqnarray} \label{S}
\deg x & \le & 1/k(\deg x +\deg y+\deg z - 1)\nonumber\\
\deg y & \le & 1/l(\deg x +\deg y+\deg z - 1)\\
\deg z & \le & 1/m(\deg x +\deg y+\deg z - 1)\,.\nonumber
\end{eqnarray}
Summing up the three inequalities in (\ref{S}), 
in virtue of our assumption, 
we obtain:
$$1/k+1/l+1/m\le (1/k+1/l+1/m-1)(\deg x +\deg y
+\deg z)\le 0\,,$$
a contradiction.

\no (b) Every one of the Platonic surfaces 
$S=S_{2,2,m},\,\, S_{2,3,3},\,\,S_{2,3,4}$ 
or $S_{2,3,5}$ admits a finite morphism 
$\C^2\to S$ (the orbit morphism of
the standard linear action on $\C^2$ of 
the corresponding finite subgroup
$\Gamma\subset {\rm SU}(2)$; see e.g., 
\cite[\S 4 and Remark 2.1]{Mil};
see also \cite[Ch. I]{Schw, Kl} or 
\cite[p.56]{Beu, BaDw} for explicit formulas). 
Since the affine plane $\C^2$ is 
$\A_1$-rich so is $S$ 
(see Remark 3 preceding the lemma).
This proves the lemma. \qed  

\medskip The next lemma is a simple 
corollary of Theorem 1 in \cite{Sch} 
(see also \cite{Br, Ve, FlZa} for relevant results).

\medskip \no {\bf Schmidt's Lemma 
\cite{Sch}.} 
{\it Let $S$ be a surface in 
$\C^3$ given by the equation 
$$z^m=f_d(x,y)$$ where $f_d\in\C[x,y]$ 
is a homogeneous polynomial of degree 
$d$ without multiple roots. Suppose that 
$m\ge 4$ and $d\ge 3$, or $m=3$ and 
$d\ge 5$, or $m=2$ and $d\ge 17$.
Then the surface $S$ is $\A_1$-poor 
and, moreover, every polynomial curve
$f\,:\,\C\to S$ passes through the 
singular point ${\overline 0}\in S$.}

\medskip The purpose of the next lemma 
is to strengthen 
the lemmas of Halphen and Schmidt
(cf. the examples below). Recall that a regular action 
of the multiplicative group $\C^*$ on an affine variety 
$X$ is called {\em good} if it has a unique fixed point
(called {\em vertex}),
and this fixed point is {\em elliptic}, that is, 
it belongs to the closure of any orbit. 
Let $S$ be a normal affine surface with a good $\C^*$-action;
denote $S^*=S\setminus V_0$ where $V_0$ is the vertex, and set 
$\G=S^*/\C^*$. 
If the curve $\G$ is rational then the 
singularity of the surface $S$ at the 
origin is called {\it quasirational} 
\cite{Ab} (cf. also examples in \cite{Ore2}). 

\medskip \no {\bf Lemma 7.}
{\em Let $S$ be a normal affine surface with a good $\C^*$-action.
Suppose 
that the singularity of the surface $S$ 
at the vertex $V_0\in S$ is not quasirational. 
Then any rational curve $r\,:\,\C\to S$, 
as well as any holomorphic entire curve 
$h\,:\,\C\to S$  
in the surface $S$ is contained in an orbit closure
$\overline {\C^*V}$ for a certain point $V\in S^*$. 
Consequently, 
any polynomial curve $f\,:\,\C\to S$
passes through the vertex $V_0\in S$, 
and so the surface $S$ is $\A_1$-poor.}

\medskip\no {\it Proof.} 
Let ${\widetilde \G}\to \G$
be a normalization. The rational mapping 
$g\,:\,\C \stackrel{r}{\longrightarrow}
S\stackrel{\pi'}{\longrightarrow} \G$ 
can be lifted to a morphism 
${\widetilde g}\,:\,\C\to {\widetilde \G}$, 
which is constant
because (by our assumption) the geometric genus 
$g({\widetilde \G})\ge 1$. Thus, the image
$r(\C)\subset S$ is contained in the closure 
${\overline {\cal O}}$
of an orbit ${\cal O}=\C^*V$ 
of the $\C^*$-action, as stated. 
The proof for an entire curve $h\,:\,\C\to S$
is similar. \qed

\medskip\no {\bf Remark.}  
This lemma (with the same proof) remains 
true also for meromorphic 
curves $m\,:\,\C\to S$ assuming that 
$e(\G)<0$ i.e., 
$g({\widetilde \G})\ge 2$ 
(cf. e.g., \cite{Ja, Grs}).

\medskip The following 
facts will be useful in order 
to provide examples of 
quasihomogeneous surfaces in $\C^3$
which satisfy the assumption of Lemma 7. 
For integers $a_1,\ldots,a_n$ 
denote $[a_1,\dots,a_n]={\rm lcm}\,(a_1,\dots,a_n)$ 
and
$(a_1,\dots,a_n)={\rm gcd}\,(a_1,\dots,a_n)$, whereas 
$\overline{(a_1,\dots,a_n)}$ denotes the vector
with the coordinates $a_1,\ldots,a_n$. 

\medskip \no {\bf Lemma 8.} 
{\em Let $f\in\C[x,y,z]$ 
be a quasihomogeneous 
polynomial such that
$$f(\lambda^{q_0}x,\,\lambda^{q_1}y,\,\lambda^{q_2}z)=
\lambda^{d}f(x,y,z)\qquad\forall\lambda\in\C$$
where $(q_0,q_1,q_2)=1$ and $q_0,\,q_1,\,q_2,\,d>0$.  
Suppose that $d\equiv 0\mod q_i, \quad i=0,1,2$, and that
the surface $S:=f^{-1}(0)\subset\C^3$ has an 
isolated singularity at the 
origin. Then the singularity $(S,\,{\bar 0})$ 
is quasirational if and only if one of the 
following two conditions holds:
\begin{enumerate}
\item[(i)] $d=[q_0,q_1,q_2]$, and for some natural numbers
$p,q,r,s$ coprime in pairs we have (up to a reordering):  
$\overline{(q_0,q_1,q_2)}=\overline{(pq,\,pr,\,qrs)}$. 
\item[(ii)] $d=2[q_0,q_1,q_2]$, and for some natural numbers
$p,q,r$ coprime in pairs we have:
$\overline{(q_0,q_1,q_2)}=\overline{(pq,\,pr,\,qr)}$.
\end{enumerate}}

\medskip\no {\em Proof}. By \cite[Prop.3]{OrlWa}, $\G=S^*/\C^*$ 
is a smooth curve of genus 
$$g(\G)={1\over 2}\biggl({d^2\over q_0q_1q_2}-d\biggl(
{1\over [q_0,q_1]}+{1\over [q_0,q_2]}+{1\over [q_1,q_2]}
\biggl)  +2    \biggl)\,.$$ 
Thus $g(\G)=0$ if and only if
\be\label{genus} 
{1\over [q_0,q_1]}+{1\over [q_0,q_2]}+{1\over [q_1,q_2]}=
{d\over q_0q_1q_2}+{2\over d}\,.
\ee
Letting $q_{ij}:=(q_i,q_j)$ we can write 
\be\label{numbers} q_0=q_{01}q_{02}q_0',\quad
q_1=q_{01}q_{12}q_1',\quad
q_2=q_{02}q_{12}q_2'\ee
where the integers 
$q_{01},\,q_{02},\,q_{12},\,q_0',\,q_1',\,q_2'$
are coprime in pairs, 
since by our assumption $(q_0,q_1,q_2)=1$. 
Set $d_0:=q_{01}q_{02}q_{12}$; then
$[q_0,q_1,q_2]=d_0q_0'q_1'q_2'$, and so 
$d=\rho d_0q_0'q_1'q_2'$ with $\rho\in\N$ 
and $q_0q_1q_2=d_0^2q_0'q_1'q_2'$.
Therefore, (\ref{genus}) can be written as
$$\rho^2-\biggl(
{1\over q_0'q_1'}+{1\over q_0'q_2'}+{1\over q_1'q_2'}
\biggl)\rho+{2\over q_0'q_1'q_2'}=0\,.$$
It follows that $\rho=1$ or $\rho=2$. 
In the first case the only solutions 
of this Diophantine equation 
are (up to a reordering) 
$\overline{(q_0',\,q_1',\,q_2')}=\overline{(1,1,s)}\quad (s\in\N)$,
and in the second case 
$\overline{(q_0',\,q_1',\,q_2')}=\overline{(1,1,1)}$ 
is the only solution. Letting in (\ref{numbers})
$q_{01}=:p,\,q_{02}=:q,\,q_{12}=:r$
and taking into account the above observations,
we get the desired conclusion. 
\qed

\medskip\no The next corollary in the 
particular case of the Pham-Brieskorn 
surfaces can be found in 
\cite[p.117]{BarKa}\footnote{Cf. \cite[Thm.2]{Ev} 
where 
in fact, several possibilities covered 
by the condition (ii$'$) below have been omited, 
because of a gap in the reduction of problem B to
problem C.}.

\medskip\no {\bf Corollary.} {\em Assume 
that the polynomial
$$f=ax^k+by^l+cz^m+\ldots\in\C[x,y,z]
\qquad (\mbox{where}\quad a,b,c\in\C^*)\,$$
is quasihomogeneous and such that the surface
$S:=f^{-1}(0)\subset\C^3$ has an 
isolated singularity at the 
origin. Then the singularity $(S,\,{\bar 0})$ 
is quasirational if and only if one of the 
following two conditions holds:
\begin{enumerate}
\item[(i$'$)] up to a reordering, $(k,lm)=1$, or  
\item[(ii$'$)] 
$(k,l)=(k,m)=(l,m)=2$. 
\end{enumerate}}

\medskip\no {\em Proof}. Letting 
$k=\rho k',\,l=\rho l',\,m=\rho m'$
where $\rho:=(k,l,m)$, set 
$d_0':=(k',l')(k',m')(l',m')$ and
$$q_0:={l'm'\over d_0'},\quad q_1:=
{k'm'\over d_0'},\quad q_2:=
{k'l'\over d_0'}\,.$$
These are the unique weights with $(q_0,q_1,q_2)=1$
making $f$ a quasihomogeneous polynomial 
of degree
$$d:=kq_0=lq_1=mq_2=\rho {k'l'm'\over d_0'}=
\rho [k',l',m']=[k,l,m]\,.$$ 
To apply Lemma 8 assume that 
$$\overline{(q_0,q_1,q_2)}=\overline{\biggl({l'm'\over d_0'},\,
{k'm'\over d_0'},\,
{k'l'\over d_0'}\biggl)}=\overline{(pq,\,pr,\,qrs)}$$
with $p,\,q,\,r,\,s$ coprime in pairs. Then
we would have 
$[q_0,q_1,q_2]=pqrs=[k',l',m']$, 
whence
$d=\rho [q_0,q_1,q_2]$. 
In view of this observation, 
it is easily seen that the condition 
(i) resp., (ii) of Lemma 8 is fulfilled 
if and only if
(i$'$) resp., (ii$'$) holds
(more precisely, iff $\rho=1$ and 
up to a reordering, $\overline{(k,l,m)}=\overline{(rs,qs,p)}$
resp., $\rho=2$ and 
up to a reordering, $\overline{(k,l,m)}=2\overline{(r,q,p)}$). 
Now the statement follows from Lemma 8. \qed 

\medskip\no {\bf Examples. 1.} 
By Lemma 7 and the above corollary, 
the Fermat cubic surface 
$x^3+y^3+z^3=0$ in $\C^3$ is $\A_1$-poor, 
and moreover, 
any rational curve in it has the diagonal 
form 
$t\longmapsto (\varphi (t) x_0,\,\varphi(t) y_0,\,\varphi (t)z_0)$
with $\varphi \in\C(t)$. 
In contrast, the cubic surface 
$x^3+y^3+z^3=1$ in $\C^3$ is rich with rational curves, 
as it is rational \cite{Man}. 
Furthermore, being non-rational \cite{ClGr},
the affine Fermat cubic threefold 
$x^3+y^3+z^3+u^3=0$ in $\C^4$ 
is unirational, and even is dominated by 
the affine space $\C^3$; 
see \cite[\S 3]{PaVa} for explicit 
formulas\footnote{In a discussion with the authors 
H. Flenner conjectured 
that this threefold does not
admit a nontrivial $\C_+$-action;
however, so far we do not possess 
a proof of this.}.   

\smallskip\no {\bf 2.} 
For the Pham-Brieskorn surfaces 
$S_{k,l,m}$, Lemma 7 and the above 
corollary provide
information additional to those given 
by Halphen's Lemma.
Namely, such a surface possesses a 
non-diagonal polynomial curve (i.e., not of the form 
$$t\longmapsto (\varphi_0^{q_0}(t),\,
\varphi_1^{q_1}(t),\,\varphi_2^{q_2}(t))\qquad {\rm with}\quad \varphi_i\in\C[t],\,\,i=0,1,2{\rm)}$$ 
if and only if one of the conditions
(i$'$) or (ii$'$) holds. (Cf. \cite{BoMu, Beu, DarGr, Ev} 
for similar results, 
including the more general situation of 
Pham-Brieskorn type surfaces 
over function fields. Besides, in \cite{DarGr} 
one can find a historical account on the subject.) 

\smallskip\no {\bf 3.} 
Let a surface $S=\{z^m-f_d(x,y)=0\}$ 
in $\C^3$ be as in Schmidt's Lemma. 
We may choose a coordinate system in 
$\C^2_{x,y}$ in such a way that 
neither $x$ nor $y$ divides the polynomial 
$f_d$. 
Then the assumptions of the above corollary are 
fulfilled with $k=l=d$, and so 
the singularity $(S,\,\bar 0)$ 
is quasirational if and only if
either $d=2$ or $(m,d)=1$. According to Lemma 7, 
the conclusion of Schmidt's Lemma remains 
true for any pair $(m,d)$ with
$m\ge 2,\,\,d\ge 2$, except for possibly 
the pairs $(2,2),\,(3,2),\,(3,4)$ and 
$(2,2k+1),\,\,k=1,\dots,7$. 

\medskip\no {\em Added in proof.} After this paper was written, 
it was established in \cite{FlZa} that 
a normal affine surface $S$ with a good $\C^*$-action 
which possesses a closed rational curve  
not passing through the vertex $V_0\in S$, has at most 
rational singularity
$(S,\,V_0)$ at the vertex. In particular, this nicely fits
Halphen's Lemma; indeed, the singularity 
$(S_{k,l,m},\,\bar 0)$ of 
the Pham-Brieskorn surface $S_{k,l,m}$
is rational precisely for the Platonic surfaces. 
This also implies that the conclusion of
Schmidt's Lemma is true
exactly when $d\ge 3$ and $(d,m)\neq (3,2)$.

\end{document}